\documentclass[twoside,draft,reqno]{birkart}

\usepackage{amssymb}

\newcommand{\bo}[1]{\mathbf{#1}}
\newcommand{\lra}{\longrightarrow}
\newcommand{\espe}[1]{{\mathrm Sp}^\infty{#1}}

\newcommand{\ra}{\rightarrow}
\newcommand{\orb}[1]{{\bo O}_{#1}}
\newcommand{\orbop}[1]{{\bo O}_{#1}^{\mathrm op}}

\newcommand{\spaces}{\bo {Spaces}}

\newcommand{\fr}{{\bo{FMod}}_{\mathrm R}}
\newcommand{\si}{\Sigma}

\newcommand{\fz}{{\bo {FMod}}_{\mathrm Z}}

\newtheorem{theorem}{Theorem}[section]
\newtheorem{lemma}[theorem]{Lemma}
\newtheorem{corollary}[theorem]{Corollary}
\newtheorem{proposition}[theorem]{Proposition}

\newtheorem{remark}[theorem]{Remark}

\newtheorem{definition}[theorem]{Definition}

\newtheorem{notation}[theorem]{Notation}

\begin{document}

\setcounter{page}{1}

\title[Recognition Principle for GEMs]
{Recognition Principle for\\
Generalized Eilenberg--Mac~Lane Spaces}

\author[B.~Badzioch]{Bernard Badzioch}

\address{
Department of Mathematics,\br
University Notre Dame,\br
Notre Dame, IN 46556\br
USA}

\email{badzioch.1@nd.edu}

\begin{abstract}
We give a homotopy theoretical characterization
of  generalized Eilenberg--Mac~Lane spaces which 
resembles the $ \Gamma $-space structure used by 
Segal to describe infinite loop spaces.
\end{abstract}

\maketitle

\section{Introduction}

A generalized Eilenberg--Mac~Lane space (GEM) is a 
space weakly equivalent to a product of  
Eilenberg--Mac~Lane spaces $\prod_{i}^{\infty}
{K(\pi_i, i)}$ with $ \pi_i $ abelian. The goal of 
this note is to prove the following characterization
of GEMs:

\begin{theorem}
\label{main}
Let $\mathrm R$ be a commutative ring with a unit and 
let ${\fr}$ be the category of finitely 
generated free $ \mathrm R$-modules. If 
$H\colon {\fr}\ra \spaces$
is a functor such that
\begin{itemize} 
\item $H(0)$ is a contractible space,
\item$ H(R)$ is connected and for every $ n$ the projections
$ {\mathrm pr}_{k}\colon {\mathrm R}^{n}\ra {\mathrm R}$
induce a weak equivalence
$ H({\mathrm R}^n) \stackrel{\simeq}{\lra} 
H({\mathrm R})^n$ 
\end{itemize}
then the space
$ H({\mathrm R})$ is weakly equivalent to a product
$ \prod_{i=1}^{\infty}
{{\mathrm K}({\mathrm M}_i , i)}$ where $ {\mathrm M}_i$
is an $ {\mathrm R}$-module.
\end{theorem}
\begin{notation}
By $ \spaces$ above and in the rest of this paper
we denote the category of simplicial sets. Consequently,
by  'space' we always mean an object of this category.
\end{notation}

The above description of GEMs is modeled after the 
$ \Gamma $-space structure introduced by Segal 
in \cite{Segal} to characterize infinite loop spaces.
Just as for infinite loop spaces one  gets
the following corollary which is implicitly 
present in the work of Bousfield~\cite{Bou} and 
Dror~\cite{Dror-cell} who apply it to 
localization functors. 
 
\begin{corollary}
If $F\colon\spaces\ra\spaces$ is a functor preserving
weak equivalences and preserving 
products up to weak equivalence then $F$ preserves GEMs. 
\end{corollary}

The rest of the paper is organized as follows. 
In section 2 the Grothendieck construction
on a diagram of small categories is recalled
and  some of its properties are stated.
In section 3 we give a description suitable
for our purposes of $\espe{X} $ -- the infinite symmetric
product on a space $ X$. In section 4 the proof of 
theorem~\ref{main} is presented. It essentially 
amounts to showing that for a functor 
$ H\colon\fr\ra\spaces$ satisfying the assumptions
of the theorem the space
$H(R)$ is a homotopy retract of $ \espe H(R)$.
Since $\espe H(R)$ is a GEM and the class of 
GEMs is closed under homotopy
retractions the claim of the theorem will follow. 

I am indebted to W. G. Dwyer for suggesting the
topic of this work to me and for many discussions. The 
referees comments were very helpful in improving the
presentation of the paper. I also want to thank 
Sylwia Zab{\l}ocka for her support and friendship.

\section{The Grothendieck construction}

\begin{definition}
Let $ \bo{Cat}$ denote the category of small 
categories and let $ \bo N$ be the telescope category:
$${\bo N}:= (0\stackrel{i_0}{\lra}1\stackrel{i_1}{\lra} \dots
\stackrel{i_{n-1}}{\lra}n\stackrel{i_n}{\lra}\dots)$$
For a functor $P:{\bf N}\rightarrow \bo{Cat}$ the 
Grothendieck construction on $ P$~\cite{Tho} is 
the category $ Gr(P)$ whose objects are pairs $ (n,c)$
where $ n \in \bo N$ and $ c\in P(n)$. A morphism 
$ (n,c)\ra (n',c')$ is a pair $(i, \varphi)$ where
$ i$ is the unique morphism $ n\ra n'$ in $ \bo N$ 
and $ \varphi \in Mor_{P(n')}(P(i)c, c')$. The 
composition of $ (i, \varphi)\colon(n,c)\ra(n',c')$
and $ (i', \varphi')\colon(n',c')\ra(n'',c'')$ is 
defined to be the pair 
$ (i'\circ i, \varphi'\circ P(i')\varphi)
\colon(n,c)\ra(n'',c'')$.
\end{definition}

For every $ n\in \bo N$ there is a functor 
$$ I_n\colon P(n)\ra Gr(P)\ \ ,\ \ I_n(c):= (n,c) $$
which lets us identify $ P(n)$ with a subcategory 
of $ Gr(P)$. It follows that any functor 
$ F\colon Gr(P)\ra \bo C$ defines a sequence of functors
$ F_n\colon P(n)\ra \bo C$, $ n=0,1,\dots$ Moreover,
for $ n\in \bo N$ and $ c\in P(n)$ let $ \beta_{n,c}$
be the image under $ F$ of the morphism 
$ (i_n, {\rm id}_{P(i_n)(c)})\in 
Mor_{Gr(P)}((n,c), (n+1, P(i_n)(c)))$. It is easy to check 
that the morphisms $\{\beta_{n,c}\}_{c\in P(n)}$ define a natural transformation of functors
$$\beta_n\colon F_n \ra F_{n+1}\circ P(i_{n}) $$
The converse is also true~\cite[A.9]{Ch-S}: any sequence 
of functors $\{F_n\colon P(n)\ra {\bo C\}}_{n\geq 0}$
and natural transformations 
$\{\beta_n\colon F_n \ra F_{n+1}\circ P(i_n)\}_{n\geq 0}$
can be used to define a functor 
$ F\colon Gr(P)\ra \bo C$, such that $ F\mid_{P(n)}= F_n.$

\begin{proposition}
\label{thomason}
For any functor $ F\colon Gr(P)\ra \bo C$ the natural 
morphism 
$$ colim_{\bo N}colim_{P(n)}F_{n}{\ra}colim_{Gr(P)}F$$
is an isomorphism.
Moreover, if $ {\bo C} = \spaces$ then the natural 
map 
$$ hocolim_{\bo N}hocolim_{P(n)}F_{n}{\ra} hocolim_{Gr(P)}F$$
is a weak equivalence. 
\end{proposition}

\begin{proof}
The first statement follows directly from the 
definition of $ Gr(P)$. The proof of the second can 
be found in~\cite[Prop. 0.2]{Sl} 
or~\cite[Cor. 24.6]{Ch-S}. 
\end{proof}

\section{Infinite symmetric products}

Let $ \si_n$ be the permutation group of the set 
$\{1,\dots, n\}$. We will denote by $ \orb{\si_n}$ 
the orbit category of $ \si_n$
whose objects are sets $ \si_n/G$ for 
$ G\subseteq\si_n$ and whose morphisms are $ \si_n$-- equivariant maps $ \si_n/G \ra \si_n/H$. Let 
$ \orbop{\si_n}$ be the opposite category. 
We can identify $ \si_n$ with the subgroup of
all elements of $ \si_{n+1}$ which leave the 
element $ n+1$ fixed. The inclusion $ \si_n\subset\si_{n+1}$ induces a functor
$$ J_{n}:\orbop{\Sigma_{n}}\ra
\orbop{\Sigma_{n+1}}\ \ ,\ \ J_n(\Sigma_{n}/G):= \Sigma_{n+1}/G$$
This data in turn can be used to define 
a functor $ O\colon {\bo N}\ra {\bo {Cat}}$:
$$O(n):= \orbop{\si_{n}}\ \ ,\ \ O(i_{n}):=J_{n}$$
Let $ \spaces_\ast$ denote the category of pointed spaces
and let $ X\in \spaces_\ast$. The group $ \si_n$ acts
on $ X^n$ by permuting the coordinates. As usual 
we have the fixed point functor
$$ F_nX\colon\orbop{\si_n}\ra \spaces_\ast$$
defined by $F_nX(\si_n/G):=(X^n)^G$ -- the fixed point
set of the action of $G$ on~$X^n$. 

\begin{remark}
\label{fixspace}
For $ G \subseteq \si_n$ let 
$|orb G|$ denote 
the number of orbits of the action of $ G$ on $\{1,\dots,n\}$. Then there is a natural isomorphism
$ (X^n)^G\cong X^{|orb G|}$.
\end{remark}

Using the embedding $ \si_n\subset\si_{n+1}$ one 
can think of $ G\subseteq\si_n$ as a subgroup
of $ \si_{n+1}$. There is an obvious isomorphism
$$ (X^{n+1})^G\cong(X^n)^G\times X$$
and since $X$ is a pointed space, 
we have a map 
$$ \beta_{n,G}\colon (X^n)^G\ra(X^{n+1})^G$$ 
$$(x_1,\dots,x_n)\mapsto(x_1,\dots,x_n,\ast)$$
One can check that the morphisms 
$ \{\beta_{n,G}\}_{G\subseteq\si_n}$ give a natural 
transformation 
$$ \beta_n\colon F_nX\ra F_{n+1}X\circ J_n$$
and from the remarks made in section 2 it follows that 
the functors $ \{F_n\}_{n\geq 0}$ and natural 
transformations $ \{\beta_n\}_{n\geq 0}$
can be assembled to define a functor 
$$ FX\colon Gr(O)\ra\spaces_\ast$$

\begin{lemma}
\label{espe}
$hocolim_{Gr(O)}FX\simeq \espe X$
\end{lemma}

\begin{proof}
By proposition \ref{thomason} we have a weak equivalence
$$ hocolim_{Gr(O)}FX\simeq hocolim_{\bo N}
hocolim_{\orbop{\si_n}}F_n X$$
Moreover, by \cite[Ch.4, Lemma A.3]{Dror-cell}
$$ hocolim_{\orbop{\si_n}}F_nX
\simeq colim_{\orbop{\si_n}}F_nX$$
(this follows from the fact that 
$F_nX$ is a free
diagram of spaces 
\cite[2.7]{Dror-diag}, and that for free diagrams
homotopy colimits coincide with ordinary colimits).
But $ colim_{\orbop{\si_n}}F_nX
\cong{\mathrm Sp}^nX$ and so we have
$$ hocolim_{\bo N}hocolim_{\orbop{\si_n}}F_nX\simeq
hocolim_{\bo N}{\mathrm Sp}^nX\simeq 
colim_{\bo N}{\mathrm Sp}^nX\cong \espe{X} $$
where the second equivalence is a consequence of
\cite[Ch. XII, 3.5]{BK}. 
 
\end{proof}

\section{The proof of theorem \ref{main}}

Let $ H\colon\fr\ra\spaces$ be a functor as in the 
theorem. One can assume that $ H$ takes its values in 
the category $ \spaces_\ast$ of pointed spaces (if not,
replace $H$ with $\widetilde{H}$, $\widetilde{H}(M) = cofib(H(0)
\ra H(M))$ for $M\in\fr$). Moreover, once we know that
the theorem holds for $R= \mathrm Z$, the ring of 
integers, the embedding $ \mathrm Z\hookrightarrow R$
will induce a functor $\fz\ra\fr$, and so the space 
$ H(R)$ will have the structure of a GEM. 
Furthermore the action of the ring $ R$ on its free module 
$ R\in \fr$ via multiplications will 
induce an action of $ R$  on $ H(R)$ and so
the homotopy groups $ \pi_i(H(R))$ will be 
$ R$--modules as claimed. 
Therefore for the rest of this paper we will assume 
that $ R = \mathrm Z$ and that 
$ H\colon\fz\ra\spaces_\ast$.

For a free abelian group on $ n$ generators
$ {\mathrm Z}^n\in\fz$ the group $ \si_n$ acts
on $ {\mathrm Z}^n$ by permuting the set of generators.
For $G\subseteq\si_n $ let $ ({\mathrm Z}^n)^G$
be the subgroup of all elements of $ {\mathrm Z}^n$
which are fixed by the action of $ G$. 

\begin{remark}
\label{fixgroup}
It is not difficult to check that, using the notation of
\ref{fixspace}, there is a natural isomorphism of groups
$ ({\mathrm Z}^n)^G\cong {\mathrm Z}^{|orb G|}$.
\end{remark}

For any $ n\in N $ we have a functor
$$ F_n{\mathrm Z}\colon \orbop{\si_n}\ra\fz\ \ ,
\ \ F_n{\mathrm Z}(\si_n/G):=({\mathrm Z}^n)^G$$
Arguments similar to those in section 3
show that one can define a functor 
$$F{\mathrm Z}\colon Gr(O)\ra\fz $$
such that 
$ F{\mathrm Z}\mid_{\orbop{\si_n}} = 
F_n{\mathrm Z}$.
Observe that $ \si_1 $ is a trivial group and so 
$ F{\mathrm Z}\mid_{\orbop{\si_1}}=F_1{\mathrm Z}$
is the constant functor with value $ \mathrm Z$.

\begin{lemma}
\label{transformation}
Let $ \bar{ Z}\colon Gr(O)\ra\fz$ be the constant
functor with value $ \mathrm Z$. There exists a 
natural transformation 
$$ \theta\colon F{\mathrm Z}\ra\bar{ Z}$$
such that $ \theta\mid_{\orbop{\si_1}}$ is an isomorphism.
\end{lemma}

\begin{proof}
For $ \si_n/G\in \orbop{\si_n}\subset Gr(O)$ define
$$\theta_{\si_n/G} \colon 
({\mathrm Z}^n)^G\ra {\mathrm Z}$$
$$(k_1,k_2,\dots,k_n)\mapsto \sum_{i}{k_i} $$
It is easy to check that these maps give the 
required transformation of functors. 
\end{proof}

Let $ H\colon \fz\ra \spaces_\ast$ be a functor
satisfying the conditions of theorem~\ref{main};
that is the projections $ {\mathrm Z}^n\ra {\mathrm Z}$
induce weak equivalences 
$ H({\mathrm Z}^n)\ra H({\mathrm Z})^n$.

\begin{lemma}
\label{espehz}
$hocolim_{Gr(O)}H\circ F{\mathrm Z}
\simeq \espe H({\mathrm Z})$ 
\end{lemma} 

\begin{proof}
It follows from \ref{fixspace} and 
\ref{fixgroup}
that for $ G\subseteq\si_n$ we have isomorphisms
$$ H(({\mathrm Z}^n)^G)\cong H({\mathrm Z}^{|orb G|})$$
and
$$ (H({\mathrm Z})^n)^G\cong H({\mathrm Z})^{|orb G|}$$
Their composition with the map 
$H(\mathrm Z^{|orb G|})\ra H(\mathrm Z)^{|orb G|}$
induced by projections $ Z^{|orb G|}\ra Z$ gives a map
$ \varphi_{n,G}\colon H(({\mathrm Z}^n)^G) 
\ra (H({\mathrm Z})^n)^G$ which, in view of the 
properties of $ H$, is a weak equivalence.
Moreover, the maps $ \{\varphi_{n,G} \}_{n\geq 0, G\subseteq\si_n}$
define a natural transformation of functors
$$ \varphi\colon H\circ F{\mathrm Z}
\ra FH(\mathrm Z)$$
Therefore we have a weak equivalence
$$hocolim_{Gr(O)}H\circ F{\mathrm Z}
\stackrel{\simeq}{\lra}
hocolim_{Gr(O)}FH({\mathrm Z}) $$
But by lemma \ref{espe} 
$hocolim_{Gr(O)}FH({\mathrm Z})\simeq
\espe H({\mathrm Z})$.
\end{proof}

To conclude the proof of the theorem observe
that the natural transformation ~$ \theta$ 
from lemma \ref{transformation} gives a 
transformation 
$$ H(\theta)\colon H\circ F{\mathrm Z}
\ra H\circ\bar{ Z}$$ 
and so induces a map 
$$hocolim_{Gr(O)}H\circ F{\mathrm Z}\ra
hocolim_{Gr(O)}H\circ\bar{ Z}\ra 
colim_{Gr(O)}H\circ\bar{ Z}\cong H({\mathrm Z})$$
On the other hand the inclusion 
$ \orbop{\si_1}\subseteq Gr(O)$ gives a map
$$ H({\mathrm Z}) \simeq 
hocolim_{\orbop{\si_1}}H\circ F_1{\mathrm Z}\ra 
hocolim_{Gr(O)}H\circ F{\mathrm Z}$$
and since $ \theta\mid_{\orbop{\si_1}}$ is an 
isomorphism it is easy to see that the composition
$$ H({\mathrm Z})\ra 
hocolim_{Gr(O)}H\circ F {\mathrm Z}
\ra H({\mathrm Z })$$
has to be a weak equivalence. But by \ref{espehz}
$hocolim_{Gr(O)}H\circ F{\mathrm Z}\simeq 
\espe H({\mathrm Z})$, and so $H({\mathrm Z})$ 
must be a GEM as a homotopy retract of a GEM 
(see ~\cite[Ch.4, Thm B.2]{Dror-cell}).
\begin{remark}
\label{final}
The above proof remains valid  
if we replace $\fz$ by the category of
free, finitely generated abelian monoids.
\end{remark}

\end{document}